\documentclass[letterpaper, 10 pt, conference]{ieeeconf}  

\IEEEoverridecommandlockouts                              

\usepackage{graphics} 
\usepackage{epsfig} 
\usepackage{mathptmx} 
\usepackage{times} 
\usepackage{amsmath} 
\usepackage{amssymb} 
\usepackage{xcolor}

\newtheorem{thm}{Theorem}
\newtheorem{remark}{Remark}

\title{\LARGE \bf Occupation Measures in the ODE Method: Wasserstein Convergence and Rates

}

\author{Dimitrios Katselis
\thanks{Dimitrios Katselis is with the ECE and CSL, Grainger College of Engineering,
        University of Illinois Urbana-Champaign, 61801, Urbana, IL, USA 
        {\tt\small katselis@illinois.edu}}
}

\begin{document}

\maketitle
\thispagestyle{empty}
\pagestyle{empty}

\begin{abstract}
We study occupation measures associated with stochastic approximation recursions within the ordinary differential equation (ODE) method. Under Lyapunov conditions ensuring that the limiting ODE has a globally exponentially stable equilibrium, we establish that the continuous-time occupation measures associated with the piecewise linear interpolation of the iterates converge almost surely in the $2$-Wasserstein topology to the Dirac mass at the equilibrium,  the unique invariant probability measure of the limiting semiflow. The choice of the $2$-Wasserstein distance  is naturally aligned with our continuous- and discrete-time Lyapunov drift analysis. We also consider two discrete-time sequences of occupation measures, induced by the  algorithmic time scale and the iteration count, respectively, and we establish their almost sure convergence to the same invariant measure in the  $2$-Wasserstein distance. We then derive corresponding  convergence rates in expectation, providing a nonasymptotic complement to the almost sure results. These expectation bounds also imply some tail probability estimates. Finally, we verify that our results apply to a canonical class of linear recursions. 
\end{abstract}

\section{INTRODUCTION}

Stochastic approximation (SA) algorithms correspond to a fundamental class of recursive schemes for solving fixed-point and optimization problems using noisy samples.  Such algorithms arise naturally when exact function evaluations or gradients are not available and only noisy or sample-based information is provided,  as in adaptive control and reinforcement learning. An important approach for the analysis of convergence and stability of SA algorithms is the ordinary differential equation (ODE) method, which interprets the continuous-time interpolation of the SA iterates as tracking trajectories of the limiting ODE generated by the mean field. Equivalently, the ODE method considers the SA updates as a noisy Euler discretization of the aforementioned limiting ODE  \cite{b96, b21, ge24, ky03}. In \cite{b96, b21, ky03}, SA algorithms are studied using a combination of weak convergence techniques, Lyapunov stability methods, invariant-set analysis, and, in generalized settings with set-valued mean fields, differential inclusions. Moreover, under standard assumptions, the limit sets of SA recursions are almost surely connected internally chain transitive invariant sets of the limiting semiflow, i.e., each such set is invariant under the semiflow, compact and connected, and any two of its points can be connected by an $(\epsilon, T)$-chain whose chain points lie entirely within the set, for every $\epsilon>0$ and $T>0$.  Analyses of recursive stochastic algorithms and connections to adaptive control and system identification are also presented in \cite{l77, lpw92}.  

More recent work has focused on the nonasymptotic analysis of SA algorithms. Due to space limitations, we present here a selection of relevant papers to illustrate  this framework. Some of these papers adapt martingale concentration methods to stochastic approximation by employing contractivity and stability properties to control deviations.  \cite{cb26} develops concentration bounds for one-step temporal difference (TD) learning or TD$(0)$,  with function approximation, by relying on Poisson-equation techniques to handle Markovian noise and relaxed concentration inequalities to address potentially unbounded iterates. Concentration results for contractive stochastic approximation under subgaussian additive or bounded multiplicative noise are presented in \cite{cmz25}. For linear stochastic approximation and TD  learning, \cite{sy19} derives finite-time bounds on error moments using Lyapunov drift arguments, including bounds on the mean square deviation from the equilibrium of the associated linear ODE. In \cite{hlr19, hk14, rss12}, bounds for stochastic gradient methods in strongly convex stochastic optimization are obtained. \cite{ls18} considers constant stepsize linear SA with iterate averaging and provides instance-dependent mean square error bounds under stability and moment conditions. High probability bounds for two-timescale linear SA are provided in \cite{dtsm18},  under martingale difference noise with iterate-dependent bounds and vanishing stepsizes. \cite{brs21}  relates TD learning to stochastic gradient descent and provides mean square error bounds. An approximation to Least Squares TD solution is considered in \cite{pkm21}  and finite-time bounds in expectation and probability are derived. Finally, \cite{tb19} derives a concentration bound for nonlinear SA using Alekseev's formula to relate the interpolated SA trajectory to the associated limiting ODE. 

An aspect that has received relatively less attention  is the study of occupation measures, whose asymptotic behavior is associated with the underlying ergodic behavior. In the ODE method,  under standard assumptions, weak limit points of occupation measures are invariant probability measures of the limiting semiflow. This allows the study of time averages by employing results from the ergodic theory of dynamical systems.  An important contribution is \cite{bs00}, which establishes ergodic properties of weak asymptotic pseudotrajectories for semiflows. In this framework, for processes that asymptotically track the semiflow of a dynamical system in the weak asymptotic pseudotrajectory sense, almost surely every weak$^*$ limit point of the associated occupation measures   is an invariant probability measure of the semiflow. Thus, under almost sure tightness of the SA occupation measures, the same conclusion applies if the interpolated SA trajectory is a weak asymptotic pseudotrajectory of the limiting semiflow.  

 Motivated by \cite{bs00} and the related discussion in \cite{b21}, we employ the ODE method to study  Wasserstein convergence and nonasymptotic rates for SA occupation measures in continuous and discrete time, instead of pathwise deviations of the iterates. The continuous-time occupation measures are defined using the piecewise linear interpolation of the SA iterates. In discrete time, we consider occupation measures induced by the algorithmic time scale and  by discrete-time averaging. When the limiting ODE has a globally exponentially stable equilibrium $x^{*}$, the limiting semiflow has a unique invariant probability measure which is atomic, namely the Dirac mass $\delta_{x^{*}}$. In this case, under standard  assumptions, we establish the almost sure convergence of these occupation measures in the $2$-Wasserstein distance to $\delta_{x^{*}}$.  Squared $W_2$ errors for $\delta_{x^{*}}$ are averaged squared deviations, directly linking the Wasserstein formulation with our Lyapunov analysis.   Finally, in the same case,  we derive  rates in expectation for their $2$-Wasserstein convergence, which also imply some tail probability bounds. Implications for the $1$-Wasserstein distance are also provided.

  \section{Preliminaries}
  
  On a filtered probability space $(\Omega, \mathcal{F}, \{\mathcal{F}_n\}_{n\geq 0}, \mathbb{P})$, we consider SA recursions of the form
 \begin{equation}\label{eq:SA1}
 x_{n+1}=x_{n}+a(n)(h(x_n)+M_{n+1}), \ \ n\geq 0,
 \end{equation}
 where $(x_n)_{n\geq 0}$ is an $\mathbb{R}^d$-valued adapted process, $h:\mathbb{R}^d \rightarrow \mathbb{R}^d$  is a Borel measurable function, $(a(n))_{n\geq 0}$ is a sequence of deterministic stepsizes, and $(M_{n})_{n\geq 1}$ is an $\mathbb{R}^d$-valued noise sequence adapted to the filtration $\{\mathcal{F}_n \}_{n\geq 0}$. Here,  $\mathcal{F}_n=\sigma(x_0,M_1,\ldots, M_n), n\geq 0$ \cite{b21, ky03}. We further  introduce the following standard assumptions \cite{b21}: 
 
   \begin{enumerate}
 \item[(A1)] (\textbf{Stepsizes})  The sequence of stepsizes $(a(n))_{n\geq 0}$ satisfies the Robbins-Monro conditions
 \begin{equation*}
 a(n)>0, \forall n\geq 0, \ \ \sum_{n\geq 0} a(n)=\infty, \ \ \sum_{n\geq 0} a(n)^2<\infty.
 \end{equation*}
 \item[(A2)] (\textbf{Regularity of the mean field $h$})  $h:\mathbb{R}^d\rightarrow \mathbb{R}^d$ is globally Lipschitz, i.e., there exists $0<L<\infty$ such that $\|h(x)-h(y)\|\leq L \|x-y\|, \forall x,y\in \mathbb{R}^d$. Then, the autonomous limiting ODE associated with (\ref{eq:SA1}), $\dot{x}(t)=h(x(t))$, is well-posed, i.e., for every initial condition $x(0)=\tilde{x}\in \mathbb{R}^d$, it admits a unique global solution $x(\cdot)\in C^1([0,\infty); \mathbb{R}^d)$ that depends continuously on $\tilde{x}$. This follows from the Picard-Lindel\"{o}f theorem \cite{b21, t12}. The associated semiflow is denoted by $(\Phi_t)_{t\geq 0}$.
  \item[(A3)] (\textbf{Martingale difference noise})   $(M_{n})_{n\geq 1}$ is a martingale difference sequence with respect to $(\mathcal{F}_n)_{n\geq 0}$, i.e., in addition to being adapted, $M_n$ is integrable for every $n\geq 1$, and   $\mathbb{E}[M_{n+1}| \mathcal{F}_n]=0$ a.s. for every $n\geq 0$.
 
 \item[(A4)] (\textbf{Conditional quadratic growth condition on the noise})  $M_n\in L^2(\Omega; \mathbb{R}^d)$ for every $n\geq 1$, and there exists a constant $0<K<\infty$ such that 
 \begin{equation*}
 \mathbb{E}[\|M_{n+1}\|^2 |\mathcal{F}_n]\leq K (1+\|x_n\|^2)\ \ \mathrm{a.s.}, \ \ n\geq 0.
 \end{equation*}
 \end{enumerate}
 We do not assume almost sure boundedness of the iterates a priori. 
 Boundedness will be established later by appropriate Lyapunov conditions on the mean field $h$. 
  
 Define the algorithmic time scale 
 $t(0)=0$  and $t(n)=\sum_{k=0}^{n-1}a(k)$, $n\geq 1$.
Let  $\bar{x}:[0,\infty)\rightarrow \mathbb{R}^d$ denote the continuous-time, piecewise linear interpolation of the SA iterates, defined by $\bar{x}(t(n))=x_n$ and 
 \begin{equation}\label{eq:SA2}
 \bar{x}(t)=x_n+ \frac{t-t(n)}{a(n)} (x_{n+1}-x_n)
 \end{equation}
  for $t\in [t(n), t(n+1)), n\geq 0$. $\bar{x}(\cdot)$ has tracking properties relative to solutions of the limiting ODE uniformly over every finite time horizon, a.s. on $\{\sup_{n\geq 0} \|x_n\|<\infty\}$  \cite{b96, b21}.  
  
  Moreover, let $\mathcal{P}(\mathbb{R}^d)$ denote the set of Borel probability measures on $\mathbb{R}^d$ endowed with the topology of weak convergence. Let also $C_0(\mathbb{R}^d)$ be the Banach space of continuous functions on $\mathbb{R}^d$ vanishing at infinity, equipped with the supremum norm. $\mathcal{P}(\mathbb{R}^d)$ is canonically embedded as a convex subset of $C_0(\mathbb{R}^d)^{*}$, the continuous dual space of $C_0(\mathbb{R}^d)$, via $\rho\rightarrow \left(f\rightarrow \int_{\mathbb{R}^d}fd\rho\right)$ for $\rho\in \mathcal{P}(\mathbb{R}^d)$. 
  By the Riesz representation theorem,  $C_0(\mathbb{R}^d)^{*}$ is isometrically isomorphic  to the space of finite signed Radon measures on $\mathbb{R}^d$ equipped with the total variation norm \cite{f99}. We associate to $\bar{x}(\cdot)$ the family of occupation measures $(\nu_t)_{t\geq 0}$ defined by 
 \begin{equation*}
 \int f d\nu_t=\frac{1}{t}\int_0^t f(\bar{x}(s)) ds, \ \ \forall f\in C_0(\mathbb{R}^d),\ \ t>0,
 \end{equation*} 
 or equivalently,
 $\nu_t=\frac{1}{t}\int_0^t \delta_{\bar{x}(s)} ds, t>0$, with  $\nu_0=\delta_{\bar{x}(0)}$. Here,  $\delta_{x}$ denotes the Dirac measure at $x\in \mathbb{R}^d$. By construction, each $\nu_t$ is a probability measure on $\mathbb{R}^d$, with $\nu_t(\mathbb{R}^d)=1$ for all $t\geq 0$. Moreover, $\nu_t$ is well-defined since $\bar{x}(\cdot)$ is measurable by noting that $\bar{x}(\cdot) \in C([0,\infty); \mathbb{R}^d)$. Under suitable conditions ensuring the almost sure boundedness of the SA iterates, $\bar{x}(\cdot)$ takes values in a path-dependent compact subset of $\mathbb{R}^d$ almost surely, and hence, each $\nu_t$ is supported on the same set. Therefore, under such conditions, almost surely $(\nu_t)_{t\geq 0}$ is tight and thus, relatively compact in $\mathcal{P}(\mathbb{R}^d)$ by Prokhorov's theorem \cite{b21}. Also, $\rho\in \mathcal{P}(\mathbb{R}^d)$ is invariant under $(\Phi_t)_{t\geq 0}$ if 
 \begin{equation*}
 \int_{\mathbb{R}^d} f(x) \rho(dx)=\int_{\mathbb{R}^d} f(\Phi_t(x)) \rho(dx), \ \ \forall t\geq 0, \ \ \forall f\in C_b(\mathbb{R}^d),
 \end{equation*}
  or equivalently, if $(\Phi_t)_{\#}\rho=\rho, \forall t\geq 0$. Here, $C_b(\mathbb{R}^d)$ is the space of bounded continuous functions on $\mathbb{R}^d$ and $(\Phi_t)_{\#}\rho=\rho \circ \Phi_t^{-1}$  is the pushforward of $\rho$ under $\Phi_t$.  Invariant probability measures and their ergodic decompositions characterize the ergodic behavior of the limiting ODE.  If $\sup_{n\geq 0}\|x_n\|<\infty$ a.s., then, almost surely, every weak limit point  $\nu_{*}$ of $(\nu_t)_{t\geq 0}$ as $t\rightarrow \infty$ is invariant under $(\Phi_t)_{t\geq 0}$ \cite{b21}. On any path where $(\nu_t)_{t\geq 0}$ is tight, its weak and weak$^{*}$ limit points coincide.
  
  Recall that the limiting ODE associated with the SA recursion (\ref{eq:SA1}) is 
  \begin{equation}\label{eq:limitingODE}
  \dot{x}(t)=h(x(t)).
  \end{equation} 
  If (\ref{eq:limitingODE}) has a globally exponentially stable equilibrium $x^{*}$, then  the unique probability measure invariant under $(\Phi_t)_{t\geq 0}$ is $\delta_{x^{*}}$. 
  
  In this paper, we consider $(\nu_t)_{t\geq 0}$ and also discrete-time occupation measures, specifically the stepsize weighted occupation measures, induced by the algorithmic time scale and given by
 $\tilde{\mu}_n=\frac{1}{t(n)}\sum_{k=0}^{n-1} a(k) \delta_{x_k}$,  $n\geq 1$,  and
 the uniformly weighted occupation measures, induced by discrete-time averaging and given by
 $\mu_n=\frac{1}{n}\sum_{k=0}^{n-1} \delta_{x_k}$, $n\geq 1$.  To study convergence of probability measures we use Wasserstein distances. For $p\geq 1$, let $\mathcal{P}_p(\mathbb{R}^d)=\{\rho \in \mathcal{P}(\mathbb{R}^d): \int_{\mathbb{R}^d} \|x\|^p \rho(dx)<\infty\}$. The $p$-Wasserstein distance between $\mu,\nu \in \mathcal{P}_p(\mathbb{R}^d)$ is defined by 
 $W_p(\mu,\nu)=\inf_{\gamma \in \mathcal{C}(\mu,\nu)}\left(\int_{\mathbb{R}^d\times \mathbb{R}^d} \|x-y\|^p \gamma(dx, dy)\right)^{1/p}$,
 where $\mathcal{C}(\mu,\nu)$ is the set of all  couplings of $\mu$ and $\nu$ \cite{v09}. In the following, our main focus is $p=2$.  Some implications for $p=1$ are also provided.

 \section{Asymptotic results}
 \label{sec:AsymptoticResults} 
 
 We present convergence results for SA occupation measures under Lyapunov stability conditions on (\ref{eq:limitingODE}). 
 
  \begin{thm}\label{thm:Asymptotic}
 Assume that (A1)-(A4) hold. Suppose that there exist $x^{*}\in \mathbb{R}^d$, a function $V:\mathbb{R}^d\rightarrow [0,\infty)$ with $V\in C^2(\mathbb{R}^d)$, and constants $c_i>0, i=1,\ldots, 5$ such that
 \begin{enumerate}
 \item[(a)] $h(x^{*})=0$
 \item[(b)] $c_1 \|x-x^{*}\|^2\leq V(x)\leq c_2 \|x-x^{*}\|^2, \ \ \forall x \in \mathbb{R}^d$
 \item[(c)] $\|\nabla V(x)\|\leq c_3 (1+\|x\|), \ \  \|\nabla^2 V(x)\|\leq c_4, \ \ \forall x \in \mathbb{R}^d$
 \item[(d)] $\langle \nabla V(x), h(x) \rangle\leq - c_5 V(x), \ \ \forall x\in \mathbb{R}^d$
 \end{enumerate}
 where $\|\cdot\|$ is the Euclidean norm for vectors or the operator norm for matrices. Assume also that $\mathbb{E}[V(x_0)]<\infty$. Then $x^{*}$ is a globally exponentially stable equilibrium of (\ref{eq:limitingODE}), $x_n\rightarrow x^{*}$ a.s. as $n\rightarrow \infty$, and $\sup_{n\geq 0} \|x_n\|<\infty$ a.s. Moreover,  $\nu_t \in \mathcal{P}_2(\mathbb{R}^d)$, $\forall t\geq 0$, $\tilde{\mu}_n, \mu_n \in \mathcal{P}_2(\mathbb{R}^d)$, $\forall n\geq 1$,  $\delta_{x^{*}}\in \mathcal{P}_2(\mathbb{R}^d)$, and additionally,  $W_2(\nu_t, \delta_{x^{*}})\rightarrow 0$ a.s. as $t\rightarrow \infty$,  $W_2(\tilde{\mu}_n,\delta_{x^{*}})\rightarrow 0$ a.s. as $n\rightarrow \infty$ and  $W_2(\mu_n,\delta_{x^{*}})\rightarrow 0$ a.s. as  $n\rightarrow \infty$.
 \end{thm}
 
 The proof is provided in Section \ref{sec:ProofThm1}.

 \begin{remark} 
 In $(c)$, the gradient bound follows from the Hessian bound, which is explicitly assumed rather than inferred from global exponential stability \cite{v23}.
 \end{remark} 
 \begin{remark}
 Since $\mathcal{P}_2(\mathbb{R}^d)\subset \mathcal{P}_1(\mathbb{R}^d)$ and $W_1(\mu,\nu)\leq W_2(\mu,\nu)$, $\forall \mu, \nu \in \mathcal{P}_2(\mathbb{R}^d)$, Theorem \ref{thm:Asymptotic} also implies that $W_1(\nu_t, \delta_{x^{*}})\rightarrow 0$ a.s. as $t\rightarrow \infty$, $W_1(\tilde{\mu}_n, \delta_{x^{*}})\rightarrow 0$ a.s. as $n\rightarrow \infty$, and $W_1(\mu_n,\delta_{x^{*}})\rightarrow 0$ a.s. as $n\rightarrow \infty$. 
 \end{remark}
  
 \section{Convergence Rates}
 \label{sec:RateResults}
 
 We  now provide rates in expectation. We consider the additional stepsize regularity condition
 \begin{equation}\label{eq:Conv1}
 \sum_{k=n_0}^{n-1} a(k)^2 e^{-\lambda (t(n)-t(k+1))}\leq B a(n), \ \ n\geq n_0
 \end{equation}
 for some constants $\lambda>0, B>0$ and an integer $n_0\geq 0$. This condition controls the exponentially weighted accumulation of squared stepsizes relative to the current stepsize.

 \begin{thm}\label{thm:Expectations}
  Assume that (A1)-(A4) hold and that  $x^{*}$, $V:\mathbb{R}^d\rightarrow [0,\infty)$ with $V\in C^2(\mathbb{R}^d)$, and $c_i>0, i=1,\ldots, 5$, satisfy conditions $(a)$-$(d)$ in Theorem \ref{thm:Asymptotic} with $\mathbb{E}[V(x_0)]<\infty$.
 Furthermore, suppose that condition  (\ref{eq:Conv1}) holds for some $\lambda \in (0,c_5), B>0$ and an integer $n_0\geq 0$, and that there exists $\tilde{K}>0$ such that, for every  $n\geq n_0$, $(1-\lambda a(n))\geq 0$ and
 \begin{equation}\label{eq:Conv2}
 \mathbb{E}[V(x_{n+1})|\mathcal{F}_n]\leq (1-\lambda a(n)) V(x_n)+\tilde{K} a(n)^2\ \ \ \mathrm{a.s.}
 \end{equation}
 
 Then the following hold. 
 
 \begin{enumerate}
 \item[(i)] For $n\geq n_0$, 
 \begin{equation*}
 \mathbb{E}[V(x_n)]\leq e^{-\lambda (t(n)-t(n_0))} \mathbb{E}[V(x_{n_0})]+ \tilde{K} B a(n),
 \end{equation*}
 and therefore, 
 \begin{equation*}
 \mathbb{E}[\|x_n-x^{*}\|^2]\leq \frac{e^{\lambda t(n_0)} \mathbb{E}[V(x_{n_0})]}{c_1} e^{-\lambda t(n)}+ \frac{\tilde{K}B}{c_1} a(n).
 \end{equation*}
 \item[(ii)] For $t\geq 0$,  
 \begin{equation}\label{eq:Conv3}
 \mathbb{E}[\|\bar{x}(t)-x^{*}\|^2]\leq C_1 e^{-\lambda t}+ C_2 r(t),
 \end{equation}
 where  $C_1>0, C_2>0$ are some constants and $r(t)=\sup\{a(n): t(n+1)>t\}$. 
 \item[(iii)] For $t>0$,
 \begin{equation*}
 \mathbb{E}[W_2(\nu_t, \delta_{x^{*}})]\leq \left(\frac{C_1}{\lambda t}+ \frac{C_2}{t}\int_0^t r(s) ds\right)^{1/2}.
 \end{equation*}
 
 \item[(iv)] For $n\geq 1$ and $C_3=\mathbb{E}[V(x_{n_0})]+\lambda \sum_{k=0}^{n_0-1}a(k)\mathbb{E}[V(x_k)]$,
  \begin{equation*}
 \mathbb{E}[W_2(\tilde{\mu}_n, \delta_{x^{*}})]\leq \left(\frac{C_3}{\lambda c_1 t(n)}+ \frac{\tilde{K}}{\lambda c_1 t(n)} \sum_{k=0}^{n-1} a(k)^2\right)^{1/2}.
 \end{equation*}
 
 \item[(v)] For $n\geq 1$ and some constant $C_4>0$, 
   \begin{equation}\label{eq:Conv_part_v}
\mathbb{E}[W_2(\mu_n,\delta_{x^{*}})]\leq \left(\frac{C_4}{ n} \sum_{k=0}^{n-1}e^{-\lambda t(k)}+ \frac{\tilde{K}B}{c_1 n} \sum_{k=0}^{n-1}a(k)    \right)^{1/2}.
 \end{equation}

   \end{enumerate}
  In particular,  if $(a(n))_{n\geq 0}$ is nonincreasing, then $r(t)=a(n)$ for $t\in [t(n), t(n+1))$, and hence $\int_0^{\infty}r(s)ds=\sum_{n\geq 0}a(n)(t(n+1)-t(n))=\sum_{n\geq 0} a(n)^2<\infty$, i.e., $r\in L^{1}([0,\infty))$. Then,  $(iii)$ yields  $\mathbb{E}[W_2(\nu_t, \delta_{x^{*}})]\lesssim t^{-1/2}$. Also, $(iv)$ yields $\mathbb{E}[W_2(\tilde{\mu}_n, \delta_{x^{*}})]\lesssim t(n)^{-1/2}$. Finally, since $e^{-\lambda t(n)}\rightarrow 0$ and $a(n)\rightarrow 0$ as $n\rightarrow \infty$,  $\frac{1}{n} \sum_{k=0}^{n-1} e^{-\lambda t(k)}\rightarrow 0$ and $\frac{1}{n}\sum_{k=0}^{n-1} a(k)\rightarrow 0$ as $n\rightarrow \infty$, and thus $(v)$ yields $\mathbb{E}[W_2(\mu_n, \delta_{x^{*}})]\rightarrow 0$ as $n\rightarrow \infty$.  
  \end{thm} 
  
 The proof is provided in Section \ref{sec:ProofThm2}.  
 \begin{remark}
The squared $W_2$ errors are averaged squared deviations and control  fractions of time or iterations  outside neighborhoods of $x^{*}$. The ODE identifies the invariant probability measure. The proofs do not use  separate tracking arguments. We also  note that (\ref{eq:Conv2}) is related to the Lyapunov drift inequality (\ref{eq:ProofThm1_4}) in the proof of Theorem \ref{thm:Asymptotic}. The bounds in Theorem \ref{thm:Expectations} complement statistical  inference in \cite{vah25} without Gaussian approximations.
 \end{remark}
 \begin{remark} 
 Since $\mathcal{P}_2(\mathbb{R}^d)\subset \mathcal{P}_1(\mathbb{R}^d)$ and $W_1\leq W_2$, the bounds in Theorem \ref{thm:Expectations} hold with $W_2$ replaced by $W_1$. For particular stepsize schedules, direct first-moment estimates based on part $(i)$ of Theorem \ref{thm:Expectations} could also be pursued. 
 \end{remark}
      
 \begin{remark}
 In the proofs of parts $(iii)$-$(v)$ of Theorem \ref{thm:Expectations}, we first bound $\mathbb{E}[W_2^2(\cdot, \delta_{x^{*}})]$, and then we employ $\left(\mathbb{E}[W_2(\cdot, \delta_{x^{*}})]\right)^2\leq \mathbb{E}[W_2^2(\cdot, \delta_{x^{*}})]$. Therefore, Theorem \ref{thm:Expectations} also provides tail probability bounds by noting that $\mathbb{P}(W_2(\cdot, \delta_{x^{*}})\geq \varepsilon)\leq \frac{1}{\varepsilon^2}\mathbb{E}[W_2^2(\cdot, \delta_{x^{*}})]$, $\forall \varepsilon>0$. For $W_1$, $\mathbb{P}(W_1(\cdot, \delta_{x^{*}})\geq \varepsilon)\leq \frac{1}{\varepsilon} \mathbb{E}[W_1(\cdot, \delta_{x^{*}})]$, $\forall \varepsilon >0$ by Markov's inequality.
 \end{remark}
 \begin{remark}
 For $a(n)=a_0 (n+1)^{-\alpha}, a_0>0$, condition (\ref{eq:Conv1}) holds when $\frac{1}{2}<\alpha <1$, and Theorem \ref{thm:Expectations} yields expected $W_2$ bounds $O(n^{-(1-\alpha)/2})$ for $\nu_{t(n)}, \tilde{\mu}_n$, and $O(n^{-\alpha/2})$ for $\mu_n$. For $\alpha=1$ and $\lambda a_0>1$, (\ref{eq:Conv1}) holds, and these bounds become $O((\log n)^{-1/2}) $ and $O(\sqrt{\log n /n})$, respectively. At $\alpha=1$, $\lambda a_0=1$, the sum in (\ref{eq:Conv1}) is of order $ \frac{\log n}{n}$, and hence (\ref{eq:Conv1}) fails.
 \end{remark}
 \begin{remark}
 Unlike the constant stepsize linear settings of \cite{mlwbj20,sy19}, we allow nonlinear mean fields with vanishing stepsizes and martingale difference noise. Compared with \cite{cmss20}, which assumes norm-based pseudo-contraction, we use smooth Lyapunov conditions directly, with matching polynomial stepsize mean square orders where both analyses apply. The expectation bounds using (\ref{eq:Conv1}) include the initial transient and complement the almost sure rates in \cite{kv24,ly22}.
 \end{remark}
 \begin{remark}
Let $X_1,\ldots, X_n$ be i.i.d. samples from a probability measure $\mu$ on a compact subset of $\mathbb{R}^d$.  Convergence rates for the empirical measures $\hat{\mu}_n=\frac{1}{n}\sum_{k=1}^n \delta_{X_k}$ typically depend on the ambient dimension $d$, e.g., $\mathbb{E}[W_1(\hat{\mu}_n,\mu)]\lesssim n^{-1/d} $ for $d\geq 3$. This is the curse of dimensionality. In our SA setting, the rate exponents exhibit no explicit dependence on $d$, since the limiting measure is a Dirac mass and convergence is governed by the  stability of the dynamics.This  is similar to the observation that, in the i.i.d. case, Wasserstein convergence rates may depend on the intrinsic dimension of  $\mu$ instead of the ambient dimension $d$ \cite{wb19}.
 \end{remark} 
\begin{remark}
Under summable bias contributions to the Lyapunov drift and (A4) for centered noise, the almost sure proof can be adapted to biased estimators \cite{chh26,s92}. Rates require further control. Clipping \cite{vah25} may control moments but need not preserve unbiasedness or Lyapunov conditions.
 \end{remark}
 
  \section{Application: Linear Stochastic Approximation}
  
  Consider the linear SA recursion 
  $x_{n+1}=x_n+a(n)(A x_n+M_{n+1})$, $n\geq 0$,
  where $A\in \mathbb{R}^{d\times d}$ is Hurwitz (i.e., all its eigenvalues have strictly negative real parts). Suppose that (A1), (A3) and (A4) hold. We note that  $h(x)=Ax$ is globally Lipschitz with $L=\|A\|$. Hence, (A2) also holds. Let $P$ be the unique symmetric positive definite solution of the Lyapunov equation $A^T P+PA=-Q$, where $Q$ is a given symmetric positive definite matrix \cite{k02}. Consider the Lyapunov function $V(x)=x^T P x$. Then, $V\geq 0$, $V\in C^2(\mathbb{R}^d)$, and the conditions $(b)$-$(d)$ in Theorem \ref{thm:Asymptotic} hold. More specifically,  $x^{*}=0$ and  $c_5= \frac{\lambda_{\min}(Q)}{\lambda_{\max}(P)}$.   Assume also that $x_0\in L^2(\Omega; \mathbb{R}^d)$, or equivalently, $\mathbb{E}[V(x_0)]<\infty$. Then, since $h(0)=0$,  Theorem \ref{thm:Asymptotic} applies.
  
  Furthermore, suppose that (\ref{eq:Conv1}) holds for some $\lambda \in (0,c_5)$, $B>0$ and an integer $n_0\geq 0$. For the aforementioned $V$, employing (A3), (A4) and $a(n)\rightarrow 0$ verifies  (\ref{eq:Conv2})  for some $\tilde{K}>0$, with $1-\lambda a(n)\geq 0$ for sufficiently large $n$. Increasing $n_0$, if necessary, preserves (\ref{eq:Conv1}). Therefore, Theorem \ref{thm:Expectations} applies.

  \section{Proof of Theorem \ref{thm:Asymptotic}}
  \label{sec:ProofThm1}
  
  We first establish the global exponential stability of the equilibrium $x^{*}$ for the limiting ODE. Let $x(\cdot)$ be a solution of $\dot{x}(t)=h(x(t))$ with $x(0)=z_0$. By (A2), $h$ is globally Lipschitz. Thus, the ODE admits a unique global solution for every initial condition $z_0\in \mathbb{R}^d$ \cite{b21, t12}. Along the  solution,   
  $\frac{d}{dt}V(x(t))=\langle \nabla V(x(t)), h(x(t)) \rangle\leq -c_5 V(x(t)).$
  Here, we employ $(d)$ in Theorem \ref{thm:Asymptotic}.     An application of the Comparison Lemma (Lemma 3.4) in \cite{k02} to $V(x(t))$ yields
  $V(x(t))\leq e^{-c_5 t} V(z_0), \ \ t\geq 0.$
  By combining $(b)$ in Theorem \ref{thm:Asymptotic} with the last inequality, we obtain
  $c_1\|x(t)-x^{*}\|^2\leq V(x(t))\leq e^{-c_5 t} V(z_0)\leq c_2 e^{-c_5 t}\|z_0-x^{*}\|^2.$
  Hence,  $ \|x(t)-x^{*}\|\leq \sqrt{\frac{c_2}{c_1}} e^{-c_5 t/2}\|z_0-x^{*}\|$, $t\geq 0$.
Therefore, by this bound and $(a)$, we conclude that $x^{*}$ is a globally exponentially stable equilibrium of the limiting ODE. 
  
  Next, we establish a useful Lyapunov inequality. Consider (\ref{eq:SA1}) and let $\xi_{n+1}=h(x_n)+M_{n+1}$. Pathwise, with the dependence on $\omega$ suppressed, we define the $C^2$ function $g:[0,1]\rightarrow [0,\infty)$ by $g(u)=V(x_n+ua(n)\xi_{n+1})$. By applying Taylor's formula with integral remainder \cite{t12}, i.e., $g(1)=g(0)+g'(0)+\int_0^1(1-u) g''(u) du$, we have that 
  \begin{align*}
  V(x_{n+1})&=V(x_n)+ a(n) \langle \nabla V(x_n), \xi_{n+1}\rangle\\&\ \  + a(n)^2 \int_0^1 (1-u)\xi_{n+1}^T \nabla^2V(x_n+ua(n)\xi_{n+1}) \xi_{n+1}du \\&\leq 
  V(x_n)+ a(n) \langle \nabla V(x_n), \xi_{n+1}\rangle+\frac{c_4}{2}a(n)^2 \|\xi_{n+1}\|^2.
    \end{align*}  
  Here, we employ $(c)$ in Theorem \ref{thm:Asymptotic}. Also, $x_n\in L^2(\Omega; \mathbb{R}^d)$ for every  $n\geq 0$ by induction, based on $\mathbb{E}[V(x_0)]<\infty$, $(b)$, (A2), (A4) and (\ref{eq:SA1}).  Taking the conditional expectation with respect to $\mathcal{F}_n$ and using (A3) and $(d)$, we obtain
    \begin{equation}\label{eq:ProofThm1_1}
    \mathbb{E}[V(x_{n+1})| \mathcal{F}_n]\leq (1-c_5 a(n)) V(x_n)+\frac{c_4}{2}a(n)^2 \mathbb{E}[\|\xi_{n+1}\|^2|\mathcal{F}_n],    
   \end{equation}   
   which holds a.s. for every $n\geq 0$. Moreover, by (A2) and $(a)$ in Theorem \ref{thm:Asymptotic}, we have that $\|h(x_n)\|\leq L \|x_n-x^{*}\|$, and therefore $\|h(x_n)\|^2\leq L^2 \|x_n-x^{*}\|^2$. Additionally, an application of the Cauchy-Schwarz inequality (or the triangle inequality) and the elementary inequality $(a+b)^2\leq 2a^2+2b^2, \forall a,b\in \mathbb{R}$ imply that $\|\xi_{n+1}\|^2\leq 2 \|h(x_n)\|^2+ 2 \|M_{n+1}\|^2$. Hence, taking the conditional expectation and using (A4), we obtain
  \begin{equation}\label{eq:ProofThm1_2}
   \mathbb{E}[\|\xi_{n+1}\|^2|\mathcal{F}_n]\leq 2 L^2 \|x_n-x^{*}\|^2+2 K (1+\|x_n\|^2)
   \end{equation}  
  a.s. for every $n\geq 0$.    Additionally, a similar argument gives $\|x_n\|^2=\|x_n-x^{*}+x^{*}\|^2\leq 2 \|x_n-x^{*}\|^2+2 \|x^{*}\|^2$. Moreover,  $\|x_n-x^{*}\|^2\leq \frac{1}{c_1}V(x_n)$ by $(b)$ in Theorem \ref{thm:Asymptotic}. By combining the last two inequalities, we have that  
   $1+\|x_n\|^2\leq (1+2\|x^{*}\|^2)+\frac{2}{c_1} V(x_n)$. 
   By replacing this inequality in (\ref{eq:ProofThm1_2}) and by employing again $(b)$, we obtain 
   \begin{equation}\label{eq:ProofThm1_3}
   \mathbb{E}[\|\xi_{n+1}\|^2 | \mathcal{F}_n]\leq 2K(1+2\|x^{*}\|^2)+ \frac{2L^2+ 4K}{c_1} V(x_n)
   \end{equation}
   a.s. for every $n\geq 0$. 
  Replacing  (\ref{eq:ProofThm1_3}) in (\ref{eq:ProofThm1_1}) gives 
   \begin{align}\label{eq:ProofThm1_4}
   \mathbb{E}[V(x_{n+1})|\mathcal{F}_n]& \leq \left(1+\underbrace{c_4 \frac{L^2+2K}{c_1} a(n)^2}_{u_n}\right)V(x_n)\nonumber\\&\ \ -\underbrace{c_5 a(n) V(x_n)}_{\gamma_n}+\underbrace{Kc_4(1+2\|x^{*}\|^2)a(n)^2}_{v_n}\nonumber\\&= (1+u_n) V(x_n)-\gamma_n+v_n,\ \ \mathrm{a.s.},\ \ n\geq 0,
   \end{align}
   where $\sum_{n\geq 0} u_n<\infty$ and $\sum_{n\geq 0} v_n<\infty$ by (A1).  (\ref{eq:ProofThm1_4}) corresponds to the aforementioned useful Lyapunov inequality.

    Next, we establish  the boundedness of the SA iterates and the convergence $x_n\rightarrow x^{*}$ a.s. Since $\mathbb{E}[V(x_0)]<\infty$, (\ref{eq:ProofThm1_4}) implies inductively that $V(x_n)$ is integrable for every $n\geq 0$. By applying the Robbins-Siegmund almost supermartingale theorem \cite{rs71} to (\ref{eq:ProofThm1_4}), we have that $V(x_n)\rightarrow \tilde{V}$ a.s. and $\sum_{n\geq 0} a(n) V(x_n)<\infty$ a.s. Here, $\tilde{V}$ is a finite nonnegative random variable. We can show that $\tilde{V}=0$ a.s. by contradiction. Assume that $\mathbb{P}(\tilde{V}>0)>0$. Let $\Omega_0=\{V(x_n)\rightarrow \tilde{V}\}\cap \{\sum_{n\geq 0} a(n)V(x_n)<\infty\}$ and note that $\mathbb{P}(\Omega_0)=1$ and $\mathbb{P}(\{\tilde{V}>0\}\cap \Omega_0)>0$. Then, for every  $\omega\in \{\tilde{V}>0\}\cap \Omega_0$, there exist $\kappa_{\omega}>0$ and an integer $N_{\omega}\geq 0$ such that $V(x_n(\omega))\geq \kappa_{\omega}, \forall n\geq N_{\omega}$. Therefore, for every such $\omega$, $\sum_{n\geq N_{\omega}} a(n) V(x_n(\omega))\geq \kappa_{\omega} \sum_{n\geq N_{\omega}} a(n)=\infty$ by (A1). This contradicts the definition of $\Omega_0$. We conclude that $\tilde{V}=0$ a.s. Then, $(b)$ in Theorem \ref{thm:Asymptotic} yields $\|x_{n}-x^{*}\|^2\leq \frac{1}{c_1} V(x_n)\rightarrow 0$ a.s. as $n\rightarrow \infty$, i.e., $x_n\rightarrow x^{*}$ a.s. as $n\rightarrow \infty$. Furthermore, every convergent real sequence is bounded, and hence $\sup_{n\geq0} V(x_n)<\infty$ a.s. Thus, $(b)$   also implies that $\sup_{n\geq 0}\|x_n-x^{*}\|<\infty$ a.s. Since  $\|x_n\|\leq \|x_n-x^{*}\|+\|x^{*}\|$, we conclude that $\sup_{n\geq 0}\|x_n\|<\infty$ a.s.

    Having established the boundedness of the SA iterates and their a.s. convergence to the equilibrium $x^{*}$, we now prove the same properties for the piecewise linearly interpolated trajectory in (\ref{eq:SA2}). For  $t\in [t(n), t(n+1))$, we have that $\bar{x}(t)=(1-\beta_{t,n})x_n+ \beta_{t,n} x_{n+1}$ for $\beta_{t,n}=\frac{t-t(n)}{a(n)}\in [0,1)$. Therefore, the triangle inequality gives $\|\bar{x}(t)-x^{*}\|\leq (1-\beta_{t,n}) \|x_n-x^{*}\|+ \beta_{t,n} \|x_{n+1}-x^{*}\|\leq \|x_n-x^{*}\|\vee \|x_{n+1}-x^{*}\|$. Consider an $\omega$ in the almost sure event $\{\omega\in \Omega: x_n(\omega)\rightarrow x^{*}, \sup_{n\geq 0}\|x_n(\omega)-x^{*}\|<\infty\}$. Since $t(n)\rightarrow \infty$ by (A1), $n=n(t)\rightarrow \infty$ as $t\rightarrow \infty$, and thus $x_n(\omega)\rightarrow x^{*}$ and $x_{n+1}(\omega)\rightarrow x^{*}$. Hence, we obtain     
  \begin{equation}\label{eq:ProofThm1_6}
  \bar{x}(t,\omega)\rightarrow x^{*}
  \end{equation}
  as $t\rightarrow \infty$ by the last bound. Additionally,  $\|\bar{x}(t,\omega)\|\leq  \|\bar{x}(t,\omega)-x^{*}\|+ \|x^{*}\|\leq (\|x_n(\omega)-x^{*}\|\vee \|x_{n+1}(\omega)-x^{*}\|)+\|x^{*}\|$, and since $\sup_{n\geq 0}\|x_n(\omega)-x^{*}\|<\infty$, we conclude that $\sup_{t\geq 0}\|\bar{x}(t,\omega)\|<\infty$ for each such $\omega$.

 Finally, for every $\omega\in \Omega$, $\bar{x}(\cdot, \omega)$ is continuous and hence bounded on each compact interval $[0,t]$, and thus $\nu_t(\omega)\in \mathcal{P}_2(\mathbb{R}^d)$, $\forall t> 0$. Also,  $\tilde{\mu}_n(\omega), \mu_n(\omega) \in \mathcal{P}_2(\mathbb{R}^d)$, $\forall n\geq 1$ because they have finite support, and $\nu_0(\omega)=\delta_{x_0(\omega)}$ and $\delta_{x^{*}}$ clearly belong to $\mathcal{P}_2(\mathbb{R}^d)$. We  consider again an $\omega$ in the almost sure event $\{\omega\in \Omega: x_n(\omega)\rightarrow x^{*}, \sup_{n\geq 0} \|x_n(\omega)-x^{*}\|<\infty\}$. We recall that (\ref{eq:ProofThm1_6}) also holds on this event. Then, 
 $W_2^2(\nu_t(\omega), \delta_{x^{*}})=\frac{1}{t} \int_0^{t} \|\bar{x}(s,\omega)-x^{*}\|^2 ds,$ $t>0.$
 Since the integrand is bounded and tends to $0$ by the previous results, its continuous-time Ces\`{a}ro average tends to zero, i.e.,  $W_2^2(\nu_t(\omega),\delta_{x^{*}})\rightarrow 0$ as $t\rightarrow \infty$. Thus, $W_2(\nu_t,\delta_{x^{*}})\rightarrow 0$ a.s. as $t\rightarrow \infty$. For  $(\tilde{\mu}_n(\omega))_{n\geq 1}$, we have that 
 $W_2^2(\tilde{\mu}_n(\omega), \delta_{x^{*}})=\frac{1}{t(n)}\sum_{k=0}^{n-1} a(k) \|x_k(\omega)-x^{*}\|^2,$
 where $\|x_k(\omega)-x^{*}\|^2\rightarrow 0$ as $k\rightarrow \infty$ and the sequence $(\|x_k(\omega)-x^{*}\|^2)_{k\geq 0}$ is  bounded. Let $\epsilon>0$ and select an integer $N_{\epsilon,\omega}\geq 1$ such that $\|x_k(\omega)-x^{*}\|^2\leq \epsilon$ for $k\geq N_{\epsilon,\omega}$. Then, for $n\geq N_{\epsilon,\omega}$,
 \begin{align*}
 &\frac{1}{t(n)}\sum_{k=0}^{n-1} a(k) \|x_k(\omega)-x^{*}\|^2 \leq \frac{1}{t(n)}\sum_{k=0}^{N_{\epsilon,\omega}-1} a(k)\|x_k(\omega)-x^{*}\|^2\\&\ \ + \frac{\epsilon}{t(n)}\sum_{k=N_{\epsilon,\omega}}^{n-1} a(k) \leq   \frac{1}{t(n)}\sum_{k=0}^{N_{\epsilon,\omega}-1} a(k)\|x_k(\omega)-x^{*}\|^2+\epsilon. 
 \end{align*}
 By noting that $t(n)\rightarrow \infty $ as $n\rightarrow \infty$, we obtain 
 $\limsup_{n\rightarrow \infty} W_2^2(\tilde{\mu}_n(\omega), \delta_{x^{*}})\leq \epsilon$.
 Since $\epsilon>0$ is arbitrary, we conclude that $W_2(\tilde{\mu}_n(\omega),\delta_{x^{*}})\rightarrow 0$ as $n\rightarrow \infty$, and therefore, $W_2(\tilde{\mu}_n,\delta_{x^{*}})\rightarrow 0$ a.s. as $n\rightarrow \infty$. Finally, for $(\mu_n(\omega))_{n\geq 1}$,
 $W_2^2(\mu_n(\omega),\delta_{x^{*}})=\frac{1}{n} \sum_{k=0}^{n-1} \|x_k(\omega)-x^{*}\|^2.$
 Since $\|x_k(\omega)-x^{*}\|^2\rightarrow 0$  as $k\rightarrow \infty$ and $(\|x_k(\omega)-x^{*}\|^2)_{k\geq 0}$ is a bounded sequence, the last Ces\`{a}ro average converges to 0 by a similar $\epsilon$-argument. Hence, $W_2(\mu_n(\omega),\delta_{x^{*}})\rightarrow 0$ as $n\rightarrow \infty$, and thus, $W_2(\mu_n,\delta_{x^{*}})\rightarrow 0$ a.s. as $n\rightarrow \infty$.   

 \begin{remark}
 Note that the independent coupling provides upper bounds for the squared $W_2$ distances, which suffice for almost sure convergence of all three occupation measures.
 \end{remark}

 \section{Proof of Theorem \ref{thm:Expectations}}
 \label{sec:ProofThm2}
 
  We focus on proving $(i)-(ii)$ and on bounding $\mathbb{E}[W_2^2(\cdot,\delta_{x^{*}})]$. The bounds on  $\mathbb{E}[W_2(\cdot,\delta_{x^{*}})]$  follow from  Jensen's inequality, i.e., $\mathbb{E}[W_2(\cdot,\delta_{x^*})]\leq (\mathbb{E}[W_2^2(\cdot, \delta_{x^*})])^{1/2}$.

  \emph{Proof of $(i)$}: By taking expectations to both sides in (\ref{eq:Conv2}), we obtain the recursion 
 \begin{equation}\label{eq:ProofThm2_0}
 \mathbb{E}[V(x_{n+1})]\leq (1-\lambda a(n)) \mathbb{E}[V(x_n)]+ \tilde{K} a(n)^2,\ \  n\geq n_0.
 \end{equation}    
 By iterating this recursion for $n\geq n_0$, we have that 
 \begin{align}\label{eq:ProofThm2_01}
 \mathbb{E}[V(x_n)]\leq& \left(\prod_{j=n_0}^{n-1} (1-\lambda a(j))\right) \mathbb{E}[V(x_{n_0})]\nonumber\\& + \tilde{K} \sum_{k=n_0}^{n-1} a(k)^2 \prod_{j=k+1}^{n-1} (1-\lambda a(j)).
 \end{align}    
 Using the elementary inequality $e^{z}\geq 1+z, z\in \mathbb{R}$ and $1-\lambda a(j)\geq 0$ for $j\geq n_0$,  we obtain, for $n\geq m\geq n_0$, 
 \begin{equation}\label{eq:ProofThm2_02}
 \prod_{j=m}^{n-1} (1-\lambda a(j))\leq e^{-\lambda \sum_{j=m}^{n-1} a(j)}=e^{-\lambda(t(n)-t(m))}.
  \end{equation}
  Therefore, by combining (\ref{eq:ProofThm2_01}) and (\ref{eq:ProofThm2_02}), we have that 
  \begin{align}\label{eq:ProofThm2_1}
  \mathbb{E}[V(x_n)]& \leq e^{-\lambda (t(n)-t(n_0))} \mathbb{E}[V(x_{n_0})]+\tilde{K} \sum_{k=n_0}^{n-1}a(k)^2 e^{-\lambda (t(n)-t(k+1))}\nonumber\\& \leq e^{-\lambda (t(n)-t(n_0))} \mathbb{E}[V(x_{n_0})]+\tilde{K} B a(n).
  \end{align}
 In the second inequality, we employ (\ref{eq:Conv1}). Moreover, using $(b)$ in Theorem \ref{thm:Asymptotic} and (\ref{eq:ProofThm2_1}), we obtain for $n\geq n_0$ that
 \begin{align}\label{eq:ProofThm2_12}
 \mathbb{E}[\|x_n-x^{*}\|^2]& \leq \frac{1}{c_1} \mathbb{E}[V(x_n)]\nonumber\\& \leq \frac{e^{\lambda t(n_0)} \mathbb{E}[V(x_{n_0})]}{c_1} e^{-\lambda t(n)}+ \frac{\tilde{K}B}{c_1} a(n).
 \end{align}
 
  \emph{Proof of $(ii)$}: For $t\geq t(n_0)$, there exists an integer $n\geq n_0$ such that $t(n)\leq t<t(n+1)$. Therefore, $\bar{x}(t)=(1-\beta_{t,n})x_{n}+ \beta_{t,n} x_{n+1}$ with $\beta_{t,n}=\frac{t-t(n)}{a(n)}\in [0,1)$. Moreover, $x\rightarrow \|x-x^{*}\|^2$ is convex. Thus, 
 $
 \|\bar{x}(t)-x^{*}\|^2\leq (1-\beta_{t,n}) \|x_{n}-x^{*}\|^2+ \beta_{t,n} \|x_{n+1}-x^{*}\|^2$ pathwise,  and hence, 
 \begin{align*}
 &\mathbb{E}[\|\bar{x}(t)-x^{*}\|^2]\leq (1-\beta_{t,n}) \mathbb{E}[\|x_{n}-x^{*}\|^2]+ \beta_{t,n} \mathbb{E}[\|x_{n+1}-x^{*}\|^2]\\& \leq \frac{e^{\lambda t(n_0)} \mathbb{E}[V(x_{n_0})]}{c_1} e^{-\lambda t(n)}+ \frac{\tilde{K}B}{c_1} (a(n)+a(n+1)), \ \ n\geq n_0.  
 \end{align*}
 We now note that $t(n)\leq t$ and $t-a(n)<t(n)$, and therefore $e^{-\lambda t(n)}<e^{\lambda a(n)} e^{-\lambda t}\leq \tilde{B} e^{-\lambda t}$ for $\tilde{B}=e^{\lambda \sup_{n\geq 0} a(n)}<\infty$ by $\sum_{n\geq 0} a(n)^2<\infty$, which implies that $a(n)\rightarrow 0$ as $n\rightarrow \infty$ and hence, $(a(n))_{n\geq 0}$ is bounded. Using the definition of $r(t)$, we conclude that 
 \begin{align}\label{eq:ProofThm2_2}
 \mathbb{E}[\|\bar{x}(t)-x^{*}\|^2]& \leq \frac{e^{\lambda t(n_0)} \mathbb{E}[V(x_{n_0})]}{c_1} \tilde{B} e^{-\lambda t}+ 2\frac{\tilde{K}B}{c_1} r(t).  
 \end{align}  
For $t\in [t(n), t(n+1)), 0\leq n<n_0$,  again  $\mathbb{E}[\|\bar{x}(t)-x^{*}\|^2]$$\leq (1-\beta_{t,n}) \mathbb{E}[\|x_{n}-x^{*}\|^2]$$+ \beta_{t,n} \mathbb{E}[\|x_{n+1}-x^{*}\|^2]$. Thus, (\ref{eq:Conv3}) holds for $C_1=\frac{e^{\lambda t(n_0)} \mathbb{E}[V(x_{n_0})]}{c_1} \tilde{B} \vee \tilde{M}_{n_{0}}\vee 1$ and $C_2=2\frac{\tilde{K}B}{c_1}$, where 
\begin{equation*}
\tilde{M}_{n_0}=\max_{0\leq k<n_0}\left\{e^{\lambda t(k+1)}\left(\mathbb{E}[\|x_k-x^{*}\|^2]\vee \mathbb{E}[\|x_{k+1}-x^{*}\|^2]\right)\right\}.
\end{equation*}
  Here, the maximum over the empty set is zero. 
  
  \emph{Proof of $(iii)$}:  For $t>0$, we note that 
 \begin{align*}
 &\mathbb{E}[W_2^2(\nu_t, \delta_{x^{*}})]=\frac{1}{t}\int_0^t\mathbb{E}[\|\bar{x}(s)-x^{*}\|^2] ds
 \\&\leq \frac{C_1}{ t} \int_0^\infty e^{-\lambda s} ds+ \frac{C_2}{ t} \int_0^t r(s) ds=\frac{C_1}{\lambda  t} + \frac{C_2}{ t} \int_0^t r(s) ds,
   \end{align*}   
   where the extension of  (\ref{eq:ProofThm2_2}) to $t\geq 0$ is employed.
   
     \emph{Proof of $(iv)$}:    
   $W_2^2(\tilde{\mu}_n, \delta_{x^{*}})=\int_{\mathbb{R}^d} \|x-x^{*}\|^2 \tilde{\mu}_n(dx)=\frac{1}{t(n)} \sum_{k=0}^{n-1} a(k) \|x_k-x^{*}\|^2\leq \frac{1}{t(n) c_1} \sum_{k=0}^{n-1} a(k) V(x_k),$
   where in the inequality we use $(b)$ in Theorem \ref{thm:Asymptotic}. Therefore, 
   \begin{align}\label{eq:ProofThm2_3}
   \mathbb{E}[W_2^2(\tilde{\mu}_n, \delta_{x^{*}})]\leq \frac{1}{t(n) c_1} \sum_{k=0}^{n-1} a(k) \mathbb{E}[V(x_k)].   
   \end{align}
   Also, (\ref{eq:ProofThm2_0}) gives
    $\lambda a(k)\mathbb{E}[V(x_{k})]\leq \mathbb{E}[V(x_k)]-  \mathbb{E}[V(x_{k+1})]+ \tilde{K} a(k)^2$ for $k\geq n_0$
   and by summing over $k=n_0$ to $n-1$,  
   \begin{align}\label{eq:ProofThm2_4}
   \sum_{k=n_0}^{n-1} a(k) \mathbb{E}[V(x_k)]& \leq \frac{1}{\lambda} \left(\mathbb{E}[V(x_{n_0})]-\mathbb{E}[V(x_n)]+ \tilde{K} \sum_{k=n_0}^{n-1} a(k)^2\right)\nonumber\\& \leq \frac{\mathbb{E}[V(x_{n_0})]}{\lambda}+ \frac{\tilde{K}}{\lambda} \sum_{k=0}^{n-1} a(k)^2.
   \end{align}
   Let $C_3=\mathbb{E}[V(x_{n_0})]+\lambda \sum_{k=0}^{n_0-1}a(k)\mathbb{E}[V(x_k)]$. Then,   $\forall n\geq 1$, $\sum_{k=0}^{n-1} a(k)\mathbb{E}[V(x_k)]\leq \frac{C_3}{\lambda}+\frac{\tilde{K}}{\lambda}\sum_{k=0}^{n-1}a(k)^2$, and by (\ref{eq:ProofThm2_3}), 
   $\mathbb{E}[W_2^2(\tilde{\mu}_n, \delta_{x^{*}})]\leq \frac{C_3}{\lambda c_1 t(n)}+ \frac{\tilde{K}}{\lambda c_1 t(n)} \sum_{k=0}^{n-1} a(k)^2.$

     \emph{Proof of $(v)$:}  We have $W_2^2(\mu_n, \delta_{x^{*}})= \frac{1}{n}\sum_{k=0}^{n-1}\|x_k-x^{*}\|^2$. Let $C_4=\frac{e^{\lambda t(n_0)} \mathbb{E}[V(x_{n_0})]}{c_1}\vee \max_{0\leq k<n_0}\left\{e^{\lambda t(k)}\mathbb{E}[\|x_k-x^{*}\|^2]\right\}\vee 1$, where the maximum over the empty set is zero. Then, by (\ref{eq:ProofThm2_12}) and the definition of $C_4$, $
   \mathbb{E}[\|x_k-x^{*}\|^2]\leq C_4 e^{-\lambda t(k)}+\frac{\tilde{K} B}{c_1} a(k)$ for every $k\geq 0$.      
    Therefore, for every $n\geq 1$, we obtain the bound 
   $ \mathbb{E}[W_2^2(\mu_n,\delta_{x^{*}})]\leq \frac{C_4}{n}\sum_{k=0}^{n-1} e^{-\lambda t(k)}+ \frac{\tilde{K} B}{c_1 n} \sum_{k=0}^{n-1} a(k)$.

 \addtolength{\textheight}{-12cm}

\end{document}